\documentclass[amssymb,11pt, draft]{amsart}
\usepackage{amsmath}
\usepackage{amsthm}
\usepackage{amsfonts}
\newcounter{leti}

\newenvironment{lettering}{\begin{list}{\bfseries \upshape (\alph{leti})}{\usecounter{leti}}}
{\end{list}}
\newcommand{\M}{\mathrm{M}}
\newcommand{\GL}{\mathrm{GL}}
\newcommand{\U}{\mathrm{U}}
\newcommand{\Symp}{\mathrm{Sp}}
\newcommand{\G}{\mathrm{G}}
\newcommand{\stack}[2]{{#1 \atop #2}}
\newtheorem{theorem}{Theorem}
\newtheorem{lemma}[theorem]{Lemma}
\newtheorem{definition}[theorem]{Definition}

\newtheorem{cor}[theorem]{Corollary}
\newtheorem{prop}[theorem]{Proposition}

\begin{document}

\title [Power series coefficients for probabilities] {Power series
coefficients for probabilities in finite classical groups}

\author{John R. Britnell}
\address{Pembroke College\\ Cambridge\\ CB2 1RF\\ UK}
\email{J.R.Britnell@dpmms.cam.ac.uk}

\author{Jason Fulman}
\address{Department of Mathematics\\
University of Pittsburgh\\
Pittsburgh, PA 15260}
\email{fulman@math.pitt.edu}



\thanks{Fulman received funding from NSA grant H98230-05-1-0031, NSF
grant DMS-0503901, and EPSRC grant EP/D001870/1. The authors are
grateful to Peter M. Neumann for organizing the 2005 Conference on
Limiting Probabilities in Finite Classical Groups, where this research
began, and for encouraging this collaboration.}

\thanks{Version of November 14, 2005}

\begin{abstract} It is shown that a wide range of probabilities and limiting probabilities in finite classical groups have integral coefficients when expanded as a power series in $q^{-1}$. Moreover it is proved that the coefficients of the limiting probabilities in the general linear and unitary cases are equal modulo 2. 
The rate of stabilization of the finite dimensional coefficients as the dimension increases is discussed.
\end{abstract}



\maketitle

\section{Introduction} \label{intro}

Recently there has been interest in understanding the proportions
of certain types of matrices over finite fields. For example an $n
\times n$ matrix is called {\it separable} if its characteristic
polynomial has no repeated roots, {\it semisimple} if its minimal
polynomial has no repeated roots, and {\it cyclic} if its
characteristic polynomial is equal to its minimal polynomial. Let
$s_{\M(n,q)}$, $ss_{\M(n,q)}$ and $c_{\M(n,q)}$ respectively
denote the probabilities that a random $n \times n$ matrix over
$\mathbb{F}_{q}$ is separable, semisimple, or cyclic. Let
$s_{\GL(n,q)}$, $ss_{\GL(n,q)}$, $c_{\GL(n,q)}$ denote the
corresponding probabilities for a random element of $\GL(n,q)$.

Motivated by questions in computational group theory \cite{NP1}, \cite{NP4}, Neumann and Praeger \cite{NP2} proved that for $n \geq 2$, 
\[
1-\frac{1}{(q^2-1)(q-1)} < c_{\M(n,q)} < 1 - \frac{1}{q^2(q+1)},\] \[ 1
- \frac{q^2}{(q^2-1)(q-1)}-\frac{1}{2}q^{-2} - \frac{2}{3} q^{-3} <
s_{\M(n,q)} < 1-q^{-1}+q^{-2}+q^{-3}.\] 
Let
$s_{\M(\infty,q)}, ss_{\M(\infty,q)}, c_{\M(\infty,q)}, s_{\GL(\infty,q)},
ss_{\GL(\infty,q)}, c_{\GL(\infty,q)}$ be the limits of the proportions
defined above as $n \rightarrow \infty$. Using generating function
techniques, it was proved independently in \cite{F} and \cite{W2}
that \[ s_{\M(\infty,q)} = \prod_{r \geq 1} (1-q^{-r}) , \ \ \
c_{\M(\infty,q)} = (1-q^{-5}) \prod_{r \geq 3} (1-q^{-r}) \] and \[
s_{\GL(\infty,q)} = 1-q^{-1}, \ \ \ c_{\GL(\infty,q)} =
\frac{(1-q^{-5})}{(1+q^{-3})}. \] Wall \cite{W2} obtained explicit
estimates on the convergence to these limits. Concerning the
semisimple limits, it was proved in \cite{F} that \[ ss_{\M(\infty,q)} =
\prod_{r \geq 1} (1-q^{-r}) \prod_{\stack{r \ge 1}{r \equiv 0, \pm 2\ 
(\bmod \ 5)}} (1-q^{-r+1})\] and \[  ss_{\GL(\infty,q)} = \prod_{\stack{r \ge 1}
{r \equiv 0, \pm 2\ (\bmod\ 5)}} \frac{(1-q^{-r+1})}{(1-q^{-r})} .\]
From these formulas, it is clear that the limiting probabilities,
when expanded as a series in $q^{-1}$, have integer coefficients.

The above proportions have also been studied in the unitary,
symplectic, and orthogonal groups \cite{NP3}, \cite{FNP},
\cite{GLu}. Aside from the applications to computational group
theory mentioned in the previous paragraph, there are
applications to the study of derangements in group actions
\cite{FG1}, \cite{FG2} and to random generation of simple
groups \cite{FG3}. Formulas for $s_{\G(\infty,q)},
ss_{\G(\infty,q)}, c_{\G(\infty,q)}$ appear in \cite{FNP} but are
very complicated. For example $$ ss_{\U(\infty,q)} =
(1+q^{-1}) \prod_{d \textrm{ odd}} A_{q,d}(1)^{\tilde{N}(d,q)}
\prod_{d \geq 1} B_{q^2,d}(1)^{\tilde{M}(d,q)}, $$ where
$$ A_{q,d}(1)= (1-q^{-d}) \left( 1+\sum_{m \geq 1} \frac{1}{|\U(m,q^d)|} \right), $$ $$ B_{q,d}(1) = (1-q^{-d}) \left( 1+\sum_{m \geq 1} \frac{1}{|\GL(m,q^d)|} \right), $$ and $\tilde{N}(d,q)$ and $\tilde{M}(d,q)$ enumerate certain sets of polynomials---see Section \ref{integ1} for their definition.

A problem posed in \cite{FNP} was to understand the
integrality properties of the coefficients when the limits are
expanded as series in $q^{-1}$. The separable case can be
treated by adapting any of three quite different existing
methods: Wall's combinatorial approach for $\M(d,q)$ and $\GL(d,q)$ \cite{W2},
Lehrer's representation theoretic approach \cite{Le}, or the
topological approach of Lehrer and Segal \cite{LeS}. However
the semisimple case seems difficult by these approaches. In
this paper we prove a general integrality result which can
handle all of these cases.

Another result of this paper is a relation between the coefficients for limiting probabilities in the general linear and unitary groups. For instance comparing the above formulas for limiting general linear probabilities with formulas for unitary limiting probabilities in \cite{FNP} one observes that
{\setlength{\arraycolsep}{1pt}
\begin{eqnarray*}\setlength{\leftmargin}{0pt}
s_{\GL(\infty,q)} &=& 1-q^{-1}\\
s_{\U(\infty,q)} &=& 1\!-\!q^{-1}\!-\!2q^{-3}\!+\!4q^{-4}\!-\!6q^{-5}\!+\!
14q^{-6}\!-\!28q^{-7}\!+\!52q^{-8}\!-\!106q^{-9}\! + \cdots\\
c_{\GL(\infty,q)} &=& 1-q^{-3}-q^{-5}+q^{-6}+q^{-8}-q^{-9} + \cdots\\
c_{\U(\infty,q)} &=& 1-q^{-3}-q^{-5}+q^{-6}-2q^{-7} + 3q^{-8} -5q^{-9} + \cdots\\
ss_{\GL(\infty,q)} &=& 1\!-\!q^{-1}\!+\!q^{-3}\!-\!2q^{-4}\!+\!2q^{-5}\!-\!q^{-6}\!-\!q^{-7}\!+\!
3q^{-8}\!-\!4q^{-9}\! + \cdots\\
ss_{\U(\infty,q)} &=& 1\!-\!q^{-1}\!-\!q^{-3}\!+\!2q^{-4}\!-\!2q^{-5}\!+\!5q^{-6}\!-\!
9q^{-7}\! +\!11q^{-8}\! -\! 20q^{-9}\! + \cdots
\end{eqnarray*}}
These expansions suggest, and we
prove, that quite generally the coefficients in the limiting general
linear and unitary expansions are equal modulo 2. This is interesting
because in the above cases there is a {\it simple} closed formula for
the general linear limits, but no such formula is known for the
unitary limits.

We also establish an integrality result for the coefficients
in the expansion as a series in $q^{-1}$ of the probability that
an element of a fixed group $\G$ satisfies certain properties;
again the methods of \cite{W2}, \cite{Le}, \cite{LeS} do not seem
applicable at our level of generality. This result gives a
slightly different approach to a question studied in \cite{W2},
\cite{Le}, \cite{LeS}: how quickly the coefficients in the power
series expansion of $s_{\G(d,q)}$ stabilize to the coefficients in
the expansion of $s_{\G(\infty,q)}$. Results are also given for
the cyclic case (studied in \cite{W2} for Lie algebras of type $A$). 
We complement
these sharp results by giving a general approach to stabilization
results which give reasonable bounds for a wide variety of cases.

The organization of the paper is as follows. Section \ref{integ1}
proves the integrality result for limiting coefficients, and the
parity result relating the limiting coefficients in the general
linear and unitary cases. Section \ref{integ2} proves the
integrality result for the case of a fixed group. Section
\ref{stabilize} uses results from Section \ref{integ2} to discuss
the rate at which separable and cyclic coefficients stabilize to
their limits, and proves a stabilization result for more general
probabilities.

For many of out arguments we shall assume that the reader is familiar with cycle 
indices for finite matrix groups, and some of their applications. These have been developed in \cite{K} and \cite{St} and, for classical groups, \cite{F}. For further developments, see \cite{B1},\cite{B2},\cite{B3},\cite{B4}.

\section{Integrality and parity of limiting coefficients} \label{integ1}

This section has two purposes. First, it will be shown that many fixed
$q$ large dimension limiting probabilities have integral coefficients
when expanded as a power series in $q^{-1}$. The integrality result is
established for the general linear, unitary, and symplectic
groups. There is no need to state results for orthogonal groups since
(as explained in the remark after the proof of Theorem
\ref{GeneralIntegers}) the arguments of \cite{FNP} show that the
corresponding limiting probabilities are obtained from those of the
symplectic group by multiplication by easily understood
factors. Second, it will be proved that the limiting coefficients in
the general linear and unitary cases are equal modulo 2. Our principal
tools are the simple transforms in Lemmas \ref{Transform1} and
\ref{Transform2}, and the identities in Lemma \ref{Identities}, taken
from \cite{FNP}.

\begin{lemma}\label{Transform1}
Let $f(x)$ be the formal power series $1+\sum_{i\ge 1}a_{i}x^{i}$. Then there exists a unique sequence $(b_{i})$ such that $f(x)=\prod_{i\ge 1}(1-x^{i})^{b_{i}}$. The sequence $(b_{i})$ consists entirely of integers if and only if every $a_{i}$ is an integer.
\end{lemma}
\begin{proof}
We define the exponents $b_{i}$ recursively. Define $b_{1}:=-a_{1}$. Suppose that we have defined $b_{1}, \dots b_{n}$, with the finite product 
$$
P_{n}(x):=\prod_{i=1}^{n}(1-x^{i})^{b_{i}}
$$ 
being equal to $1+\sum_{i\ge 1}c_{i}x^{i}$, and that we have done this in such a way that $c_{i}=a_{i}$ when $i\le n$. Then we can force $P_{n+1}(x)$ to agree as far as the $x^{n+1}$ coefficient by defining $b_{n+1}:=c_{n+1}-a_{n+1}$. This is sufficient to prove existence and uniqueness of the exponents $b_{i}$. It is obvious that if every $b_{i}\in\mathbb{Z}$, then $f(x)$ has integer coefficients. The proof of the converse is by induction; suppose that all the $a_{i}\in\mathbb{Z}$. Then certainly $b_{1}\in\mathbb{Z}$. Suppose that $b_{1},\dots,b_{n}$ are all integers. Then $\prod_{i=1}^{n}(1-x^{i})^{b_{i}}$ has integer coefficients, and so $c_{n+1}\in\mathbb{Z}$, and hence $b_{n+1}\in\mathbb{Z}$.
\end{proof}
\begin{lemma}\label{Transform2}
Let $(a_{i})$ be a sequence of even integers, and let $f(x)$ be the formal power series $1+\sum_{i\ge 1}a_{i}x^{i}$. Then there exists a unique sequence $(b_{i})$ such that 
$$
f(x)=\prod_{i\ge 1}\left(\frac{1-x^{i}}{1+x^{i}}\right)^{b_{i}}.
$$ The sequence $(b_{i})$ consists entirely of integers if and only if every $a_i$ is an even integer.
\end{lemma}
\begin{proof}
We first note that the function $\frac{1-x^{i}}{1+x^{i}}$, when expressed as a power series in $x^{i}$, has only even coefficients after the constant coefficient, and that the coefficient of $x^{i}$ is $-2$. 
As in the proof of Lemma \ref{Transform1} we define $b_{i}$ recursively. First put $b_{1}:= -\frac{a_{1}}{2}$. Suppose that $b_{1},\dots, b_{n}$ have been defined; then define the partial product  
$$
P_{n}(x):=\prod_{i=1}^{n}\left(\frac{1-x^{i}}{1+x^{i}}\right)^{b_{i}}.
$$
Let $P_{n}(x)=1+\sum_{i\ge 1}c_{i}x^{i}$, and suppose that $b_{1},\dots, b_{n}$ have been defined in such a way that $a_{i}=c_{i}$ when $i\le n$. Then we can force $P_{n+1}(x)$ to agree with $f(x)$ as far as the $x^{n+1}$ coefficient by defining 
$b_{n+1}:=\frac{1}{2}(c_{n+1}-a_{n+1})$. Now since $P_{n}(x)$ consists of a finite product of power series with even coefficients, it follows that $c_{n+1}$ is even. Since $a_{n+1}$ is even by stipulation, it follows that $b_{n+1}$ is an integer.
\end{proof}

We define the following quantities, which count certain sets of
polynomials. In this definition, we write $\mu$ for the arithmetic
M\"{o}bius function.

\begin{definition} Let $e(q)$ be $1$ if $q$ is even, and $2$ if $q$ is odd.
\begin{lettering}
\item $$N(d,q):=\frac{1}{d}\sum_{a|d}\mu(a)(q^{\frac{d}{a}}-1)$$
\item $$\widetilde{N}(d,q):=\left\{\begin{array}{ll}
\frac{1}{d}\sum_{a|d}\mu(a)(q^{\frac{d}{a}}+1) & \textrm{if } d \textrm{  is odd,} \\
0 & \textrm{if } d \textrm{ is even.}\end{array}\right.$$
\item $$\widetilde{M}(d,q):=\frac{1}{2}\left(N(d,q^2)-\widetilde{N}(d,q)\right)$$
\item $$N^{*}(d,q):=\left\{\begin{array}{ll}
\frac{1}{d}\sum_{\stack{a|d}{a \textrm{
\emph{odd}}}}\mu(a)\left(q^{\frac{d}{2a}}+1-e(q)\right) & \textrm{if }
d \textrm{ is even,} \\ e(q) & \textrm{if } d=1, \\ 0 & \textrm{if }
d>1,\ d \textrm{ odd.} \end{array}\right.$$
\item $$M^{*}(d,q):=\frac{1}{2}\Big(N(d,q)-N^{*}(d,q)\Big)$$
\end{lettering}
\end{definition}

The following lemma brings together several identities from \cite{FNP}, namely
Lemma 1.3.10 part (b), Lemma 1.3.14 parts (a) and (d), and Lemma 1.3.17 parts (a), (c), and (e).
\begin{lemma}\label{Identities}
Let $e(q)$ be $1$ if $q$ is even, and $2$ if $q$ is odd. Suppose that \mbox{$|x|<q^{-1}$}. Then
\begin{lettering} 
\item \label{GL-} 
$$ 
\prod_{d\ge 1}(1-x^{d})^{N(d,q)} = \frac{1-qx}{1-x} 
$$ 
\item \label{U-}  
$$ 
\prod_{d \textrm{ \emph{odd}}}(1-x^{d})^{\widetilde{N}(d,q)}
\prod_{d\ge 1}(1-x^{2d})^{\widetilde{M}(d,q)} = \frac{1-qx}{1+x}
$$
\item \label{UN-+} 
$$ 
\prod_{d \textrm{ \emph{odd}}}\left(\frac{1-x^{d}}{1+x^{d}}\right)^{\widetilde{N}(d,q)} =
\frac{(1-x)(1-qx)}{(1+x)(1+qx)}
$$
\item \label{Sp-} 
$$ 
\prod_{d\ge 1}(1-x^{d})^{N^{*}(2d,q)}\prod_{d\ge 1}(1-x^{d})^{M^{*}(d,q)} =
\frac {1-qx}{(1-x)^{e(q)}}
$$
\item \label{SpN-+} 
$$ 
\prod_{d\ge 1}\left(\frac{1-x^{d}}{1+x^{d}}\right)^{N^{*}(2d,q)} =
\frac{1-qx}{(1-x)^{e(q)-1}}
$$
\item \label{Neg}
$$
\prod_{d\ge 1}(1-x^{d})^{N^{*}(2d,q)}\prod_{d\ge 1}(1+x^{d})^{M^{*}(d,q)} =
\frac {1-qx^2}{(1-x)^{e(q)-1}(1+x)^{e(q)}}
$$
\end{lettering}
\end{lemma}

We are now able to state and prove our first results on integrality of
power series coefficients in Lemmas \ref{Integers1} and
\ref{Integers2}. Recall that an infinite product $\prod_{n} (1+r_n)$
is said to converge absolutely if $\prod_{n=1}^N (1+|r_n|)$ converges,
and that $\prod_n (1+r_n)$ converges absolutely over a domain $D$ in
the complex plane if and only if $\sum_n |r_n|$ converges over $D$.
\begin{lemma}\label{Integers1}
Let $r > q^{-1}$, and let $(a_{i})$ be a series of integers such that the product
$$
P(x):=\prod_{d\ge 1}\left(1+\sum_{i\ge 1}a_{i}x^{di}\right)^{N(d,q)}
$$ 
converges absolutely whenever $|x|<r$. Then $P(q^{-1})$ has a power series expansion in $q^{-1}$ with integer coefficients. Furthermore, if we define
$$
Q(x):=\prod_{d \textrm{ \emph{odd}}}\left(1+\sum_{i\ge 1}a_{i}x^{di}\right)^{\widetilde{N}(d,q)}
\prod_{d\ge 1}\left(1+\sum_{i\ge 1}a_{i}x^{2di}\right)^{\widetilde{M}(d,q)}
$$
and
$$
R(x):=\prod_{d \ge 1}\left(1+\sum_{i\ge 1}a_{i}x^{di}\right)^{N^{*}(2d,q)}
\prod_{d\ge 1}\left(1+\sum_{i\ge 1}a_{i}x^{di}\right)^{M^{*}(d,q)}
$$
then the expansions of $Q(q^{-1})$ and $R(q^{-1})$ also have integer coefficients.
\end{lemma}
\begin{proof}
By Lemma \ref{Transform1}, there exists a unique integer series $(b_{i})$ which satisfies \mbox{$1+\sum_{i\ge 1}a_{i}x^{i} = \prod_{i\ge 1}(1-x^{i})^{b_{i}}$}. We write
\begin{eqnarray}
\nonumber P(x) & = & \prod_{d\ge 1}\prod_{i\ge 1}(1-x^{di})^{b_{i}N(d,q)} \\
\label{Pidentity} & = &\prod_{i\ge 1}\left(\;\prod_{d\ge 1}(1-x^{di})^{N(d,q)}\right)^{b_{i}}.
\end{eqnarray}
Now for a given $i$, the product $\prod_{d\ge 1}(1-x^{di})^{N(d,q)}$
converges only when $|x^{i}|<q^{-1}$.  Since $P(x)$ is absolutely
convergent when $|x|<r$, it follows that $b_{i}=0$ for any $i$ such
that $|r|^{i} > q^{-1}$. Hence every term $\prod_{d\ge
1}(1-x^{di})^{N(d,q)}$ which is present (i.e.\ has non-zero exponent $b_i$)
in the product (\ref{Pidentity}) converges when $|x|<r$. By part
\ref{GL-} of Lemma \ref{Identities}, we obtain
$$
P(x) = \prod_{i\ge 1}\left(\frac{1-qx^{i}}{1-x^{i}}\right)^{b_{i}}.
$$
In particular, this identity is valid when $x=q^{-1}$,  
which shows that the expansion of $P(q^{-1})$ in powers of $q^{-1}$ has integer coefficients. 

The products $Q(x)$ and $R(x)$ may be treated in exactly the same way, except that instead of appealing to part \ref{GL-} of Lemma \ref{Identities}, we use part \ref{U-} for $Q(x)$ and part \ref{Sp-} for $R(x)$.
Note also that absolute convergence of $P(x)$ for $|x|<r$ implies absolute
convergence of $Q(x)$ and $R(x)$ for $|x|<r$, by the criterion mentioned before
the statement of the lemma.
\end{proof}

Lemma \ref{Integers1} is sufficient to deal with limiting probabilities in general linear groups. For the other classical groups, we require the following complementary result.
\begin{lemma}\label{Integers2}
Let $r > q^{-1}$, and let $(a_{i})$ be a sequence of even integers such that the product
$$
\prod_{d\ge 1}\left(1+\sum_{i\ge 1}a_{i}x^{di}\right)^{N(d,q)}
$$  converges absolutely for $|x|< r$. Define
\begin{eqnarray*}
A(x) & := & \prod_{d \textrm{ \emph{odd}}}
\left(1+\sum_{i \geq 1} a_{i}x^{di}\right)^{\widetilde{N}(d,q)} \\
B(x) & := & \prod_{d\ge 1}\left(1+\sum_{i \geq 1} a_{i}x^{di}\right)^{N^{*}(2d,q)}. 
\end{eqnarray*}
Then the power series expansions in $q^{-1}$ of $A(q^{-1})$ and $B(q^{-1})$
have integer coefficients.
\end{lemma}
\begin{proof}
By Lemma \ref{Transform2}, we may write
$$
1+\sum_{i \geq 1} a_{i} x^{i}=\prod_{i \geq 1} \left(\frac{1-x^{i}}{1+x^{i}}\right)^{b_{i}}.
$$
Now we may proceed as in the proof of Lemma \ref{Integers1}, making use of parts \ref{UN-+} and \ref{SpN-+} of Lemma \ref{Identities} for $A(x)$ and $B(x)$ respectively.
\end{proof}

Lemma \ref{Integers2} leads us to consider the parity of the coefficients $(a_{i})$ of the power series lying within our infinite products. This approach turns out to be fruitful in terms of proving integrality of the coefficients of the expanded power series, and will also lead to the somewhat unexpected result given in Theorem \ref{GLUparity}. The following lemma is useful in this respect.
\begin{lemma}\label{Evenintegers}
Let $(a_{i})$ and $(b_{i})$ be sequences of integers. Then the power series expansion of
$$
\frac{1+\sum_{i}a_{i}x^{i}}{1+\sum_{i}b_{i}x^{i}}
$$
has even coefficients (except for the constant coefficient) if and only if $a_{i}-b_{i}$ is even for all $i$.
\end{lemma}
\begin{proof}
In the ring $\mathbb{Z}[[x]]$ of formal power series in $x$ with coefficients from $\mathbb{Z}$, let
$\langle 2x \rangle$ be the principal ideal generated by $2x$. If $g(x)$ is invertible (i.e. has constant coefficient $1$), then we observe that 
$$
f(x)-g(x)\in \langle 2x \rangle \Longleftrightarrow g^{-1}(x)\Big(f(x)-g(x)\Big)\in \langle 2x \rangle \Longleftrightarrow
\frac{f(x)}{g(x)}\in 1+ \langle 2x \rangle.
$$
This suffices to prove the lemma.
\end{proof}

We may now bring together Lemmas \ref{Integers1} and \ref{Integers2} in the following way:
\begin{lemma}\label{Integers3}
Let $r > q^{-1}$, and let $(a_{i})$ and $(b_{i})$ be sequences of integers such that 
$a_{i}-b_{i}$ is even for all $i$, and such that the products
$$
\prod_{d\ge 1}\left(1+\sum_{i\ge 1}a_{i}x^{di}\right)^{N(d,q)}
$$
and
$$
\prod_{d\ge 1}\left(1+\sum_{i\ge 1}b_{i}x^{di}\right)^{N(d,q)}
$$
are absolutely convergent for $|x|< r$.
Define 
\begin{eqnarray*}
F(x) & := & \prod_{d \textrm{ \emph{odd}}}\left(1+\sum_{i}a_{i}x^{di}\right)^{\widetilde{N}(d,q)}
\prod_{d\ge 1}\left(1+\sum_{i}b_{i}x^{2di}\right)^{\widetilde{M}(d,q)} \\
G(x) & := & \prod_{d\ge 1}\left(1+\sum_{i}a_{i}x^{di}\right)^{N^{*}(2d,q)}
\prod_{d\ge 1}\left(1+\sum_{i}b_{i}x^{di}\right)^{M^{*}(d,q)}.
\end{eqnarray*}
Then $F(q^{-1})$ and $G(q^{-1})$ have power series expansions in $q^{-1}$ with integer coefficients.
\end{lemma}
\begin{proof}
We may rewrite $F(x)$ as 
\begin{eqnarray}
\label{F(x)} F(x) & = & \prod_{d \textrm{ odd}}
\left(\frac{1+\sum_{i}a_{i}x^{di}}{1+\sum_{i}b_{i}x^{di}}\right)^{\widetilde{N}(d,q)}
\prod_{d \textrm{ odd}}\left(1+\sum_{i}b_{i}x^{di}\right)^{\widetilde{N}(d,q)}\\ & & \nonumber\times
\prod_{d\ge 1}\left(1+\sum_{i}b_{i}x^{2di}\right)^{\widetilde{M}(d,q)}. \end{eqnarray}
Now  by Lemma \ref{Evenintegers}, the quotient
$$
\frac{1+\sum_{i}a_{i}x^{di}}{1+\sum_{i}b_{i}x^{di}}
$$
expands with even coefficients except for the constant term. It
follows from Lemma \ref{Integers2} that the first of the three
infinite products in (\ref{F(x)}) expands with integer coefficients. And
by Lemma \ref{Integers1}, so do the second and third products taken
together. This proves the result for $F(x)$; the proof for $G(x)$ is
similar.
\end{proof}

From this point on we shall be working in full generality, rather than
concentrating on the particular examples of cyclic, separable or
semisimple elements. We work directly with generating functions
derived from cycle indices. For background on cycle indices of finite
classical groups, see \cite{F}. The paper \cite{St} is another useful
reference and works out examples of cycle index calculations for
$\GL(d,q)$ and $\M(d,q)$.
\begin{definition} \label{type}
Let $\Lambda$ be a (possibly infinite) set of partitions of positive integers. Let $\alpha$ be an element of a finite classical group $\G$. 
For each monic irreducible polynomial $f$ over $\mathbb{F}_{q}$, define $\lambda_{f}(\alpha)$ to be the partition whose parts are $a_{1}\dots a_{n}$, where $f^{a_{1}},\dots,f^{a_{n}}$ are the powers of $f$ amongst the elementary divisors of $\alpha$. 
\begin{enumerate}
\item If $\G \in \{ \GL, \U \}$, we say $\alpha$ is of
$\Lambda$-type if the partitions $\lambda_{f}(\alpha)$ are all either empty or in $\Lambda$.
\item If $\G= \Symp$, we say $\alpha$ is of
$\Lambda$-type if the partitions $\lambda_{f}(\alpha)$ for $f \neq z \pm 1$ are all either empty or in $\Lambda$, and $\lambda_{z\pm 1}(\alpha)$ are empty.
\end{enumerate}
\end{definition}
\begin{definition}
Define
$\Lambda_{\G(d,q)}$ to be the probability that a randomly chosen element of
$\G(d,q)$ is of $\Lambda$-type, and
$\Lambda_{\G(\infty,q)}$ to be the limiting probability as $d$
increases, the limit being taken only over even values of $d$ if
$\G=\Symp$. (The existence of this limit is explained in the
proof of Theorem \ref{GeneralIntegers}.)
\end{definition}
{\bf Remark:} The requirement in Definition \ref{type} that $\lambda_{z \pm 1}$ should 
be empty in the symplectic case is for convenience. Limiting probabilities without
this restriction are obtained by multiplying limiting probabilities
with this restriction by a factor corresponding to $z-1$ and a factor
corresponding to $z+1$. These factors are easily understood in any
particular case.
\begin{center}
\end{center}

We define quantities $C_{\GL,\lambda}(q^d)$ and $C_{\U,\lambda}(q^d)$
which appear when working with cycle indices. These quantities are
sizes of certain centralizers, but we do not need this fact and shall
define them by formulae.
\begin{definition}\label{centralizers}
Let $\lambda$ be a partition with $m_{i}$ parts of size $i$ for all
$i$, and let
$k(\lambda)= 2 \sum_{i<j} im_im_j + \sum_i (i-1) m_i^2$ . Then
\begin{eqnarray*}
C_{\GL,\lambda}(q^d) &:=& q^{k(\lambda)d}\prod_{i}|\GL(m_{i},q^{d})|\\ & = & q^{k(\lambda)d} \prod_i q^{d {m_i \choose 2}}(q^{d m_i}-1) \cdots (q^d-1)    , \\
C_{\U,\lambda}(q^d) &:=& q^{k(\lambda)d}\prod_{i}|\U(m_{i},q^{d})|\\ &  = & q^{k(\lambda)d} \prod_i q^{d {m_i \choose 2}}(q^{d m_i}-(-1)^{m_i}) \cdots (q^d+1).
\end{eqnarray*}
\end{definition}

The following lemma will be useful.
\begin{lemma}\label{stongs}
\begin{enumerate}
\item $\displaystyle \sum_{|\lambda|=n} \frac{1}{C_{\GL,\lambda}(q^d)} =
\frac{1}{q^{nd}(1-q^{-d}) \cdots (1-q^{-nd})}$.
\item $\displaystyle \sum_{|\lambda|=n} \frac{1}{C_{\U,\lambda}(q^d)} =
\frac{1}{q^{nd}(1+q^{-d}) \cdots (1-(-1)^n q^{-nd}) }$.
\end{enumerate}
\end{lemma}
\begin{proof}  The first assertion may be found in \cite{St}; it is a consequence 
of Fine and Herstein's count of nilpotent matrices \cite{FH}. For $d=1$ the
second assertion follows from the first, since \[
C_{\U,\lambda}(q) = (-1)^{|\lambda|} C_{\GL,\lambda}(-q).\] For
general $d$ replace $q$ by $q^d$. \end{proof}

We are now in a position to establish a principal result of this
section. As mentioned earlier and as explained after the proof of
Theorem \ref{GeneralIntegers}, there is no need to state results for
orthogonal groups.

\begin{theorem}\label{GeneralIntegers}
Let $\Lambda$ be a set of partitions of positive integers.  Then
$\Lambda_{\G(\infty,q)}$ may be expressed as a power series in
$q^{-1}$ whose coefficients are integers if
\mbox{$\G\in\{\GL,\U,\Symp\}$}.
\end{theorem}
\begin{proof}
If $\Lambda$ does not contain the unique partition of $1$, it is not
hard to show that $\Lambda_{\G(\infty,q)}=0$ (this also follows from
Theorem \ref{stabthe2} in the case of $GL$). We shall therefore
suppose throughout that $\Lambda$ does contain this
partition. Throughout the proof we use the notation that if $A(u)=\sum
a_n u^n$ and $B(u) = \sum b_n u^n$, then $A(u) \ll B(u)$ means that
$a_n \leq b_n$ for all $n$. We also assume familiarity with cycle
indices of finite classical groups \cite{F}.
\begin{enumerate}
\item Suppose $\G=\GL$.
In this case we may specialize the general linear group cycle index to get
$$
1+\sum_{d\ge 1}\Lambda_{\GL(d,q)}u^{d} =
\prod_{d\ge 1}\left(1+\sum_{\lambda\in\Lambda}\frac{u^{d|\lambda|}}
{C_{\GL,\lambda}(q^d)}\right)^{N(d,q)}.
$$
By applying part \ref{GL-} of Lemma \ref{Identities} with $x=q^{-1}u$, this may 
be written as $\frac{A(u)}{1-u}$ where 
\[ A(u)= (1-q^{-1}u) \prod_{d\ge 1} \left[ (1-q^{-d}u^d) \left(
1+\sum_{\lambda\in\Lambda}\frac{u^{d|\lambda|}}
{C_{\GL,\lambda}(q^d)} \right) \right]^{N(d,q)}.\] From part 1 of
Lemma \ref{stongs} and the fact that $(1) \in \Lambda$, it is not hard
to see that \begin{eqnarray*}&&  \prod_{d\ge 1} \left[ (1-q^{-d}u^d) \left(
1+\sum_{\lambda\in\Lambda}\frac{u^{d|\lambda|}}
{C_{\GL,\lambda}(q^d)} \right) \right]^{N(d,q)} \\ & \ll & \prod_{d \ge 1} \left( 1 +
\frac{u^d}{q^d(q^d-1)} + \sum_{n \geq 2} \frac{u^{nd}}{q^{nd}(1-q^{-d}) \cdots (1-q^{-nd})} \right)^{N(d,q)}\\ & \ll & \prod_{d \ge
1} \left( 1+ \frac{2u^d}{q^{2d}} + 4 \sum_{n \geq 2}
\frac{u^{nd}}{q^{nd}} \right)^{N(d,q)}. \end{eqnarray*} The last step
used the fact from \cite{NP2} that $\frac{1}{(1-q^{-1}) \cdots
(1-q^{-n})} \leq 4$ for all $n$ and $q \geq 2$. Thus
$\frac{A(u)}{1-u}$ is analytic in an open disc of radius
$q^{\frac{1}{2}}$ except for a simple pole at $u=1$. It follows that
$\Lambda_{\GL(\infty,q)}$ is equal to the residue at that
pole, which is
$$
(1-q^{-1}) \prod_{d\ge 1}
\left[ (1-q^{-d}) \left(1+\sum_{\lambda\in\Lambda}\frac{1}
{C_{\GL,\lambda}(q^d)}\right) \right]^{N(d,q)}.
$$
We can certainly find integers $(a_{i})$ such that
$$
(1-q^{-d})\left(1+\sum_{\lambda\in\Lambda}\frac{1}
{C_{\GL,\lambda}(q^d)}\right)=1+\sum_{i\ge 1}a_{i}q^{-di}.
$$ 
An argument similar to that of the previous paragraph shows that the product
$$
F(x):=\prod_{d\ge 1}\left(1+\sum_{i\ge 1}a_{i}x^{di}\right)^{N(d,q)}
$$
converges absolutely for $|x|<q^{-\frac{1}{2}}$, and hence we may appeal to Lemma \ref{Integers1} to show that $F(q^{-1})$ has integer coefficients in its expansion; since $\Lambda_{\GL(\infty,q)}=(1-q^{-1})F(q^{-1})$, this is enough to prove this case of the theorem. 
\item Suppose next that $\G=\U$. By specializing the cycle index of 
$\U(d,q)$, we obtain the identity
\begin{eqnarray*}
& & 1+\sum_{d\ge 1}\Lambda_{\U(d,q)}u^{d} \\ & = &  
\prod_{d \textrm{ odd}}\left(1+\sum_{\lambda\in\Lambda}\frac{u^{d|\lambda|}}
{C_{\U,\lambda}(q^d)}\right)^{\widetilde{N}(d,q)}\prod_{d\ge 1}\left(1+\sum_{\lambda\in\Lambda}\frac{u^{2d|\lambda|}}
{C_{\GL,\lambda}(q^{2d})}\right)^{\widetilde{M}(d,q)},
\end{eqnarray*}
which by means of part \ref{U-} of Lemma \ref{Identities} (with $x=q^{-1}u$) can be rewritten as
\begin{eqnarray*}
\lefteqn{\frac{1+q^{-1}u}{1-u}\prod_{d \textrm{ odd}}
\left[(1-q^{-d}u^{d})\left(1+\sum_{\lambda\in\Lambda}\frac{u^{d|\lambda|}}
{C_{\U,\lambda}(q^d)}\right)\right]^{\widetilde{N}(d,q)} \times} \\
& & \prod_{d\ge 1}\left[(1-q^{-2d}u^{2d})\left(1+\sum_{\lambda\in\Lambda}\frac{u^{2d|\lambda|}}
{C_{\GL,\lambda}(q^{2d})}\right)\right]^{\widetilde{M}(d,q)}.
\end{eqnarray*}
Arguing as in the $\G=\GL$ case (but using both parts of Lemma
\ref{stongs}), one sees that apart from the explicit simple pole at
$u=1$, this is analytic in an open disc of radius
$q^{\frac{1}{2}}$. Hence the value of $\Lambda_{\U(\infty,q)}$ is
equal to its residue at $u=1$. This is equal to
\begin{eqnarray*}
\lefteqn{(1+q^{-1})\prod_{d \textrm{ odd}}
\left[(1-q^{-d})\left(1+\sum_{\lambda\in\Lambda}\frac{1}
{C_{\U,\lambda}(q^d)}\right)\right]^{\widetilde{N}(d,q)} \times} \\
& & \prod_{d\ge 1}\left[(1-q^{-2d})\left(1+\sum_{\lambda\in\Lambda}\frac{1}
{C_{\GL,\lambda}(q^{2d})}\right)\right]^{\widetilde{M}(d,q)}.
\end{eqnarray*}
We can find integer sequences $(a_{i})$ and $(b_{i})$ such that for all $d$, 
\begin{eqnarray*}
(1-q^{-d})\left(1+\sum_{\lambda\in\Lambda}\frac{1}
{C_{\U,\lambda}(q^d)}\right) & = & 1+\sum_{i\ge 1}a_{i}q^{-di}, \\
(1-q^{-d})\left(1+\sum_{\lambda\in\Lambda}\frac{1}
{C_{\GL,\lambda}(q^d)}\right) & = & 1+\sum_{i\ge 1}b_{i}q^{-di}.
\end{eqnarray*}
Let us consider $C_{\U,\lambda}(q^d)$ and $C_{\GL,\lambda}(q^d)$ as
polynomials in $q$. It is clear from the definitions of these
quantities that the difference of the coefficients of these two
polynomials is even for any given power of $q$. It follows easily that
the difference of the reciprocals $C_{\U,\lambda}(q^d)^{-1}$ and
$C_{\GL,\lambda}(q^d)^{-1}$, when expanded as a power series in
$q^{-1}$, will have even coefficients, and hence that $a_{i}-b_{i}$ is
even for all $i$. Define
$$
F(x)=\prod_{d \textrm{ odd}}\left(1+\sum_{i\ge 1}
a_{i}x^{di}\right)^{\widetilde{N}(d,q)}\prod_{d\ge 1}\left(1+\sum_{i\ge 1}
b_{i}x^{2di}\right)^{\widetilde{M}(d,q)}.
$$ 
As in the case $\G=\GL$, both factors in the product $F(x)$ converge
absolutely when $|x|<q^{-\frac{1}{2}}$. We now invoke Lemma
\ref{Integers3}, which tells us that the expansion of $F(q^{-1})$ in
powers of $q^{-1}$ has integer coefficients. But then this is also
true for $\Lambda_{\U(\infty,q)}$, which is equal to
$(1+q^{-1})F(q^{-1})$.
\item Suppose $\G=\Symp$. Specializing the cycle index of $\Symp(2d,q)$ gives the identity 
\begin{eqnarray*}
\lefteqn{1+\sum_{d\ge 1}\Lambda_{\Symp(2d,q)}u^{d}  \ =} & & \\
& & \prod_{d \ge 1}\left(1+\sum_{\lambda\in\Lambda}\frac{u^{d|\lambda|}}
{C_{\U,\lambda}(q^d)}\right)^{N^*(2d,q)}\prod_{d\ge
1}\left(1+\sum_{\lambda\in\Lambda}\frac{u^{d|\lambda|}}
{C_{\GL,\lambda}(q^d)}\right)^{M^*(d,q)}. \end{eqnarray*} Using part \ref{Sp-}
of Lemma \ref{Identities} with $x=q^{-1}u$, and arguing as in the
previous cases, one deduces that $\Lambda_{\Symp}(\infty,q)$ is equal to
\begin{eqnarray*} \lefteqn{(1-q^{-1})^{e(q)} \prod_{d \geq 1} \left[ (1-q^{-d}) \left(1+
\sum_{\lambda \in \Lambda} \frac{1}{C_{\U,\lambda}(q^d)} \right)
\right]^{N^*(2d,q)} \times}\\ & & \prod_{d \geq 1} \left[ (1-q^{-d}) \left(1+
\sum_{\lambda \in \Lambda} \frac{1}{C_{\GL,\lambda}(q^d)} \right)
\right]^{M^*(d,q)}. \end{eqnarray*} The remainder of the argument is
similar to the previous cases.
\end{enumerate}
\end{proof}


{\bf Remark:} As mentioned in the introduction, the limiting
orthogonal probabilities are simple functions of the limiting
symplectic probabilities. If $(1) \not \in \Lambda$, it is not hard to
show that all limiting probabilities are 0 (for instance one could use
an argument similar to that of Theorem \ref{stabthe2}). To handle the
case $(1) \in \Lambda$, we extend an idea from \cite{FNP} for the
cases of separable, cyclic, and semisimple matrices. Suppose for
example that the dimension of the space is even and that $\lambda_{z
\pm 1}$ are empty. Then if one considers the $0$-dimensional space to
be of positive type, the sum of the cycle indices for the positive and
negative type orthogonal groups is equal to the cycle index of the
symplectic groups. Thus $ \Lambda_{O^+(\infty,q)} +
\Lambda_{O^-(\infty,q)} = \Lambda_{\Symp(\infty,q)}$. If one lets
$X(u)$ denote the difference of the cycle indices for the positive and
negative type orthogonal groups, then arguing as in \cite{FNP} (or
using Lemmas 2.6.1 and 3.7.2 of \cite{W1}) one deduces that
\[ X(u)  = \prod_{d \ge
1}\left(1+\sum_{\lambda\in\Lambda}\frac{(-1)^{|\lambda|}
u^{d|\lambda|}} {C_{\U,\lambda}(q^d)}\right)^{N^*(2d,q)}\prod_{d\ge
1}\left(1+\sum_{\lambda\in\Lambda}\frac{u^{d|\lambda|}}
{C_{\GL,\lambda}(q^d)}\right)^{M^*(d,q)}.\] It follows that $X(u)$ is
analytic in an open disc of radius $q^{1/2}$, which implies that
$\Lambda_{O^+(\infty,q)} = \Lambda_{O^-(\infty,q)}$, and hence that
both of these probabilities are equal to
$\frac{\Lambda_{\Symp(\infty,q)}}{2}$. To prove the analyticity
assertion about $X(u)$, one uses the fact that $(1) \in \Lambda$ and
part \ref{Neg} of Lemma \ref{Identities} with $x=q^{-1} u$ to write
$X(u)$ as
\begin{eqnarray*}
& & \prod_{d \ge 1}\left(1-\frac{u^d}{(q^d+1)} + \cdots
\right)^{N^*(2d,q)}\prod_{d\ge 1}\left(1 + \frac{u^d}{(q^d-1)} + \cdots
\right)^{M^*(d,q)}\\ & = & \frac{1-q^{-1}u^2}{(1-q^{-1}u)^{e(q)-1}
(1+q^{-1}u)^{e(q)}} \prod_{d \ge 1}\left( \frac{1-\frac{u^d}{(q^d+1)} +
\cdots}{(1-\frac{u^d}{q^d})} \right)^{N^*(2d,q)}\\ & & \cdot \prod_{d\ge 1}\left(
\frac{1 + \frac{u^d}{(q^d-1)} + \cdots}{(1+\frac{u^d}{q^d})}
\right)^{M^*(d,q)}. \end{eqnarray*} Then one argues as in the $\G=\GL$
case of the proof of Theorem \ref{GeneralIntegers}.
\begin{center}
\end{center}

As a corollary of Theorem \ref{GeneralIntegers}, we answer one of the questions raised in \cite{FNP}. Note that $O$ refers to an orthogonal group on an odd dimensional space, and that $O^{\pm}$ refer to orthogonal groups on an even dimensional space.
\begin{cor}
\label{answer} The coefficients of powers of $q^{-1}$ in the limiting
probabilities $s_{\G(\infty,q)}, ss_{\G(\infty,q)}, c_{\G(\infty,q)}$
are integers for $\G \in \{\GL,\U,\Symp,O,O^+,O^-\}$ except for the
cases
\begin{enumerate}
\item $s_{O^{\pm}}$ in even characteristic,
\item $c_{O^{\pm}}$ in odd or even characteristic,
\item $ss_{O^{\pm}}$ in odd or even characteristic.
\end{enumerate} For these three cases, the coefficients are half-integers.
\end{cor}
\begin{proof} This is clear from Theorem \ref{GeneralIntegers} and the formulas for limiting probabilities in \cite{FNP}. \end{proof}

The following theorem is a somewhat curious outcome of our study of
parity. This relationship is interesting because as mentioned in the
introduction, there are simple exact formulas for the limiting
proportion of regular semisimple, cyclic, and semisimple matrices in
the general linear case, but in the unitary case it is not even known
which of these limiting proportions is a rational function of $q$.
\begin{theorem}\label{GLUparity}
Let $\Lambda$ be a set of partitions of positive integers containing $(1)$. Write 
\begin{eqnarray*}
\Lambda_{\GL(\infty,q)} & = & 1+\sum_{i\ge 1}a_{i}q^{-i} \\
\Lambda_{\U(\infty,q)} & = & 1+\sum_{i\ge 1}b_{i}q^{-i}.
\end{eqnarray*}
Then $a_{i}-b_{i}$ is even for all $i$.
\end{theorem}
\begin{proof}
Let us write
\begin{eqnarray*}
(1-q^{-d})\left(1+\sum_{\lambda\in\Lambda}\frac{1}{C_{\GL,\lambda}(q^d)}\right) 
& = & 1 + \sum_{i\ge 1}v_{i}q^{-di}, \\
(1-q^{-d})\left(1+\sum_{\lambda\in\Lambda}\frac{1}{C_{\U,\lambda}(q^d)}\right) 
& = & 1 + \sum_{i\ge 1}w_{i}q^{-di}.
\end{eqnarray*}
Then as we observed in the proof of Theorem \ref{GeneralIntegers}, $v_{i}-w_{i}$ is even for all $i$, and
\begin{eqnarray*}
\Lambda_{\U(\infty,q)} & = & (1+q^{-1})\prod_{d \textrm{ odd}}
\left(\frac{1 + \sum_{i\ge 1}w_{i}q^{-di}}{1 + \sum_{i\ge 1}v_{i}q^{-di}}\right)^{\widetilde{N}(d,q)}\times \\
& & \prod_{d \textrm{ odd}}
\left(1 + \sum_{i\ge 1}v_{i}q^{-di}\right)^{\widetilde{N}(d,q)}\prod_{d\ge 1}
\left(1 + \sum_{i\ge 1}v_{i}q^{-2di}\right)^{\widetilde{M}(d,q)}.
\end{eqnarray*}
Following our usual procedure, we use Lemmas \ref{Transform1}, \ref{Transform2}, and \ref{Evenintegers} to transform this into the form
\begin{eqnarray*}
\lefteqn{(1+q^{-1})\prod_{d \textrm{ odd}}\prod_{i\ge 1}
\left(\frac{1-q^{-di}}{1+q^{-di}}\right)^{y_{i}\widetilde{N}(d,q)}\times} \\
& & \prod_{d \textrm{ odd}}\prod_{i\ge 1}(1-q^{-di})^{z_{i}\widetilde{N}(d,q)}
\prod_{d\ge 1}\prod_{i\ge 1}(1-q^{-2di})^{z_{i}\widetilde{M}(d,q)}
\end{eqnarray*}
for integer series $(y_{i})$ and $(z_{i})$. Now we use parts \ref{UN-+} and \ref{U-} of Lemma \ref{Identities} to obtain
$$
\Lambda_{\U(\infty,q)} =
 (1+q^{-1})\prod_{i\ge 1}\left[
\frac{(1-q^{-i})(1-q^{1-i})}{(1+q^{-i})(1+q^{1-i})}\right]^{y_{i}} \prod_{i\ge 1}
\left[\frac{1-q^{1-i}}{1+q^{-i}}\right]^{z_{i}}.
$$
On the other hand, by the proof of Theorem \ref{GeneralIntegers} and part \ref{GL-} of Lemma \ref{Identities}, we may write
\begin{eqnarray*}
\Lambda_{\GL(\infty,q)} & = & (1-q^{-1}) \prod_{d \geq 1} \left(1+ \sum_{i} v_i q^{-di} \right)^{N(d,q)}\\
& = & (1-q^{-1})\prod_{i\ge 1}\left(
\frac{1-q^{1-i}}{1-q^{-i}}\right)^{z_{i}}.
\end{eqnarray*}
Therefore 
$$
\frac{\Lambda_{\U(\infty,q)}}{\Lambda_{\GL(\infty,q)}}=
\frac{1+q^{-1}}{1-q^{-1}}\prod_{i\ge 1}\left[
\left(\frac{1-q^{-i}}{1+q^{-i}}\right)^{y_{i}+z_{i}}
\left(\frac{1-q^{1-i}}{1+q^{1-i}}\right)^{y_{i}}\right],
$$
the expansion of which has even coefficients (since $\frac{1-q^{-1}}{1+q^{-1}}$ does).
The theorem now follows from Lemma \ref{Evenintegers}.
\end{proof}

\section{Integrality of Finite Dimensional Coefficients} \label{integ2}

    This section proves an integrality result for the coefficients
    of probabilities in a finite classical group when expanded as
    a power series in $q^{-1}$. Here the group is fixed, so the
    result is non-asymptotic, which removes the need to deal with
    issues of convergence. For the special cases of regular
    semisimple elements in the setting of Lie algebras rather than
    Lie groups, the paper \cite{Le} gives interpretations of this
    result in terms of topology and representation theory of the
    Weyl group. The argument presented here is very much in the
    spirit of Section \ref{integ1}, but now one needs the following
    variations of Lemmas \ref{Transform1} and \ref{Transform2}
    which involve power series in two variables.
\begin{lemma} \label{qTransform1}
Let $f(u,q^{-1})$ be the formal power series $1+ \sum_{1 \leq i,j} a_{i,j} u^i q^{-j}$. Then there exists a unique sequence
$(b_{i,j})$ such that \[ f(u,q^{-1}) = \prod_{1 \leq i,j} (1- u^i
q^{-j})^{b_{i,j}}.\]  The sequence $(b_{i,j})$ consists of integers if and only
if every $a_{i,j}$ is an integer. \end{lemma}
\begin{proof} The proof is nearly identical to that of Lemma \ref{Transform1}, except that the exponents $b_{i,j}$ are defined by induction on $n:=i+j$. Thus we define 
\begin{eqnarray*}
b_{1,1} & := &  -a_{1,1},\\
b_{k,n-k} & := & -a_{k,n-k} + c_{k,n-k},
\end{eqnarray*}
where $c_{k, n-k}$ is the coefficient of $u^{k}q^{-(n-k)}$ in the expansion of
$$
\prod_{\stack{1 \leq i,j}{\stack{i<k}{j<n-k}}} (1- u^i q^{-j})^{b_{i,j}}.
$$
It is a straightforward matter to show that the $b_{i,j}$ satisfy the statements of the lemma. 
\end{proof}
\begin{lemma} \label{qTransform2} Let $(a_{i,j})$ be
a sequence of even integers. Let $f(u,q^{-1})$ be the formal power
series $1+ \sum_{1 \leq i,j} a_{i,j} u^i q^{-j}$. Then there
exists a unique sequence $(b_{i,j})$ such that \[ f(u,q^{-1}) =
\prod_{1 \leq i,j} \left( \frac{1-u^i q^{-j}}{1+u^i q^{-j}}
\right)^{b_{i,j}}.\] The sequence $(b_{i,j})$ consists entirely of
integers if and only if every $(a_{i,j})$ is an even integer. \end{lemma}
\begin{proof} The argument is a straightforward modification of Lemma \ref{Transform2}, except that the $b_{i,j}$ are defined by induction on $i+j$, as in Lemma \ref{qTransform1}.
\end{proof}

Now the main result of this section can be proved.
\begin{theorem}
\label{GeneralIntegers2} Let $\Lambda$ be a set of partitions of
positive integers. Then $\Lambda_{\G(d,q)}$ may be expressed as a
power series in $q^{-1}$ whose coefficients are integers if $\G \in
\{\GL,\U,\Symp\}$.
\end{theorem}
\begin{proof} The three cases are treated separately. For background on cycle indices of finite classical groups, see \cite{F}.
\begin{enumerate}
\item Suppose that $\G=\GL$. The general linear group cycle index gives
$$
1+\sum_{d\ge 1}\Lambda_{\GL(d,q)}u^{d} =
\prod_{d\ge 1}\left(1+\sum_{\lambda\in\Lambda}\frac{u^{d|\lambda|}}
{C_{\GL,\lambda}(q^d)}\right)^{N(d,q)}.
$$ By Lemma \ref{qTransform1} and part \ref{GL-} of Lemma
\ref{Identities}, this can be rewritten as \[ \prod_{d \geq 1} \left(
\prod_{1 \leq i,j} (1-u^{id}q^{-jd})^{b_{i,j}} \right)^{N(d,q)} =
\prod_{1 \leq i,j} \left( \frac{1-u^i q^{1-j}}{1-u^i q^{-j}} \right)^{b_{i,j}}.\] This
implies the result in the $\G=\GL$ case.
\item Suppose that $\G=\U$. Specializing the cycle index of the unitary
groups gives that \begin{eqnarray*} & & 1+\sum_{d\ge
1}\Lambda_{\U(d,q)}u^{d} \\ & = & \prod_{d \textrm{
odd}}\left(1+\sum_{\lambda \in \Lambda}\frac{u^{d|\lambda|}}
{C_{\U,\lambda}(q^d)}\right)^{\widetilde{N}(d,q)}\prod_{d\ge
1}\left(1+\sum_{\lambda \in \Lambda}\frac{u^{2d|\lambda|}}
{C_{\GL,\lambda}(q^{2d})}\right)^{\widetilde{M}(d,q)}.
\end{eqnarray*} Using the equation $C_{\U,\lambda}(q) = (-1)^{|\lambda|}
C_{\GL,\lambda}(-q)$ and Lemma \ref{qTransform1}, this can be
written as \[ \prod_{d \textrm{ odd}} \left( \prod_{1 \leq i,j} (1-
\frac{ (-1)^{i+j} u^{id}}{ q^{jd}})^{b_{i,j}} \right)^{\tilde{N}(d,q)}
\prod_{d \geq 1} \left( \prod_{1 \leq i,j} (1- \frac{u^{2id}}{
q^{2jd}})^{b_{i,j}} \right)^{\tilde{M}(d,q)}.\] By part \ref{U-} of
Lemma \ref{Identities}, this is \[ \prod_{1 \leq i,j} \left(
\frac{1-u^i q^{1-j}}{1+u^iq^{-j}} \right)^{b_{i,j}} \prod_{d \textrm{ odd}}
\left[ \prod_{1 \leq i,j} \left( \frac{(1- (-1)^{i+j} u^{id}
q^{-jd})}{(1-u^{id}q^{-jd})} \right)^{b_{i,j}}
\right]^{\tilde{N}(d,q)} .\] When one writes \[  \prod_{1 \leq i,j} \left( \frac{(1- (-1)^{i+j} u^{id}
q^{-jd})}{(1-u^{id}q^{-jd})} \right)^{b_{i,j}} = 1 + \sum_{1 \leq i,j} c_{i,j} u^i q^{-j}, \] the $c_{i,j}$ are clearly all even. Thus, applying Lemma \ref{qTransform2} and then part \ref{UN-+} of Lemma \ref{Identities} gives that  \begin{eqnarray*} & & 1+\sum_{d\ge 1}\Lambda_{\U(d,q)}u^{d} \\
& = &  \prod_{1 \leq i,j} \left(
\frac{1-u^i q^{1-j}}{1+u^iq^{-j}} \right)^{b_{i,j}} \prod_{d \textrm{ odd}} \left[ \prod_{1 \leq i,j} \left( \frac{1-u^{id} q^{-jd}}{1+u^{id} q^{-jd}} \right)^{a_{i,j}} \right]^{\tilde{N}(d,q)}\\
& = &  \prod_{1 \leq i,j} \left(
\frac{1-u^i q^{1-j}}{1+u^iq^{-j}} \right)^{b_{i,j}} \prod_{1 \leq i,j} \left( \frac{(1-u^iq^{-j})(1-u^iq^{1-j})}{(1+u^iq^{-j})(1+u^iq^{1-j})} \right)^{a_{i,j}}. \end{eqnarray*}
This implies the result for the case $\G=\U$.
\item Suppose that $\G=\Symp$. The argument is similar to the unitary case. Using Lemma \ref{qTransform1}, then part \ref{Sp-} of Lemma \ref{Identities}, followed by Lemma \ref{qTransform2} and part \ref{SpN-+} of Lemma \ref{Identities}, it follows that
\begin{eqnarray*}
& & 1+\sum_{d\ge 1}\Lambda_{\Symp(2d,q)}u^{d}\\
& = & \prod_{d \ge 1}\left(1+\sum_{\lambda\in\Lambda}\frac{u^{d|\lambda|}}
{C_{\U,\lambda}(q^d)}\right)^{N^*(2d,q)}\prod_{d\ge
1}\left(1+\sum_{\lambda\in\Lambda}\frac{u^{d|\lambda|}}
{C_{\GL,\lambda}(q^d)}\right)^{M^*(d,q)}\\
& =& \prod_{d \ge 1} \left( \prod_{1 \leq i,j} \left( 1- \frac{(-1)^{i+j} u^{id}}{q^{jd}} \right)^{b_{i,j}} \right)^{N^*(2d,q)} \left( \prod_{1 \leq i,j} (1- \frac{u^{id}}{q^{jd}})^{b_{i,j}} \right)^{M^*(d,q)}\\
& = & \prod_{1 \leq i,j} \left( \frac{1-u^iq^{1-j}}{(1-u^iq^{-j})^{e(q)}} \right)^{b_{i,j}} \prod_{d \geq 1} \left[ \prod_{1 \leq i,j} \left( \frac{1-\frac{(-1)^{i+j}u^{id}}{q^{jd}}}{1-\frac{u^{id}}{q^{jd}}} \right)^{b_{i,j}} \right]^{N^*(2d,q)}\\
& = & \prod_{1 \leq i,j} \left(
\frac{1-u^iq^{1-j}}{(1-u^iq^{-j})^{e(q)}} \right)^{b_{i,j}} \prod_{d \geq 1}
\left[ \prod_{1 \leq i,j} \left(
\frac{1-u^{id}q^{-jd}}{1+u^{id}q^{-jd}} \right)^{a_{i,j}}
\right]^{N^*(2d,q)}\\
& = & \prod_{1 \leq i,j} \left( \frac{1-u^iq^{1-j}}{(1-u^iq^{-j})^{e(q)}} \right)^{b_{i,j}} \prod_{1 \leq i,j} \left( \frac{1-u^iq^{1-j}}{(1-u^iq^{-j})^{e(q)-1}} \right)^{a_{i,j}}.
\end{eqnarray*}
This completes the proof for the case $\G=\Symp$. 
\end{enumerate}
\end{proof}

{\bf Remark:} The case of the orthogonal groups is easily understood
using the above approach. Suppose for instance that the dimension is
even and that $\lambda_{z \pm 1}$ are empty. Then considering the 0
dimensional space to be of positive type, the sum of the cycle indices
of the positive and negative type orthogonal groups is equal to the
cycle index of the symplectic groups. So by Theorem
\ref{GeneralIntegers2}, the coefficient of $q^{-j}$ in the sum of the
$O^+(2d,q)$ and $O^-(2d,q)$ probabilities is an integer. As in the
remark after Theorem \ref{GeneralIntegers}, let $X(u)$ denote the
difference of the cycle indices for $O^+(2d,q)$ and $O^-(2d,q)$.
Then as before $X(u)$ is equal to \[ \prod_{d \ge
1}\left(1+\sum_{\lambda\in\Lambda}\frac{(-1)^{|\lambda|}
u^{d|\lambda|}} {C_{\U,\lambda}(q^d)}\right)^{N^*(2d,q)}\prod_{d\ge
1}\left(1+\sum_{\lambda\in\Lambda}\frac{u^{d|\lambda|}}
{C_{\GL,\lambda}(q^d)}\right)^{M^*(d,q)}, \] which in the notation of
the $\G=\Symp$ case of Theorem \ref{GeneralIntegers2}, is equal to
\[ \prod_{d \ge 1}
\left( \prod_{1 \leq i,j} (1-(-1)^{j} u^{id}q^{-jd})^{b_{i,j}}
\right)^{N^*(2d,q)} \left( \prod_{1 \leq i,j} (1-
u^{id}q^{-jd})^{b_{i,j}} \right)^{M^*(d,q)}. \] This is slightly
different from the cycle index of the symplectic groups since the
power of $(-1)$ in the first factor is different, but the same
argument as in the $\G=\Symp$ case of Theorem \ref{GeneralIntegers2}
shows that the power series expansion in $q^{-1}$ for the coefficient
of $u^d$ in $X(u)$ has integral coefficients. Thus the coefficient of
$q^{-j}$ in the $O^+(2d,q)$ and $O^-(2d,q)$ probabilities is a
half-integer.

\section{Rate of Stabilization of Coefficients} \label{stabilize}

    This section studies the question of how large $d$ must be so that
    the coefficient of $q^{-n}$ in $\Lambda_{\G(d,q)}$ is equal to the
    coefficient of $q^{-n}$ in $\Lambda_{\G(\infty,q)}$. Subsection
    \ref{rsc} obtains sharp results in the regular semisimple and
    cyclic cases, when the group in question is $\GL$ or $\U$. For
    these groups it is well known that an element is regular
    semisimple if and only if it is separable. Then Subsection
    \ref{general} uses themes from earlier sections of this paper to
    prove a general stabilization result; while not always sharp it is
    broadly applicable.

\subsection{Stabilization for Regular Semisimple and Cyclic Probabilities} \label{rsc}

    It should be noted that the regular semisimple case has been
    studied by several authors. Lehrer \cite{Le} obtained results in
    the setting of Lie algebras rather than Lie groups, but they were
    not sharp. Sharp results in types $A,B,C$ for the Lie algebra case
    and for the $\GL$ case appear in \cite{LeS} using topological
    methods. Wall \cite{W2} uses combinatorial techniques to obtain
    sharp results for $\GL$ and its Lie algebra for the regular
    semisimple case and cyclic case. The argument presented here has
    similarities to that of Wall \cite{W2}, but seems different enough
    to record. Results are worked out for the general linear and
    unitary groups; similar methods will apply to the symplectic and
    orthogonal groups, but this is much more laborious.

    Letting $\G$ denote $\GL$ or $\U$, we use the notation that \[
    s_{\G}(u,q) := 1 + \sum_{d \geq 1} u^d s_{\G(d,q)} \] and \[
    c_{\G}(u,q) := 1 + \sum_{d \geq 1} u^d c_{\G(d,q)}. \] Here, as in
    the introduction, $s_{\G(d,q)}$ is the proportion of separable
    elements in $\G(d,q)$ and $c_{\G(d,q)}$ is the proportion of
    cyclic elements in $\G(d,q)$.

\begin{prop} \label{stab1} Let 
$F(a,s) = \frac{1}{s} \sum_{r|s,a} \mu(r) (-1)^{s/r} {s/r + a/r - 1
\choose a/r}$. Then \[ s_{\GL}(u,q) = \frac{\prod_{a \geq 0} \prod_{s
\geq 1} (1-u^s q^{1-s-a})^{F(a,s)}}{(1+\frac{u}{q-1})}. \] \end{prop}
\begin{proof} From the cycle index of general linear groups \cite{F},\cite{St} , \begin{eqnarray*}
\left( 1+\frac{u}{q-1} \right) s_{\GL}(u,q) & =  & \prod_{d \geq
1} \left( 1+\frac{u^d}{q^d-1} \right)^{\frac{1}{d} \sum_{r|d} \mu(r) q^{d/r}}\\
& = & \exp\left( \sum_{d \geq 1} \frac{1}{d} \sum_{r|d} \mu(r) q^{d/r} \cdot
\log \left( 1+\frac{u^d}{q^d-1} \right) \right)\\ & = & \exp \left( - \sum_{d \geq 1} \frac{1}{d}
\sum_{r|d} \mu(r) q^{d/r} \sum_{i \geq 1} \frac{(-1)^i u^{id}}{i
q^{id}(1-q^{-d})^i} \right). \end{eqnarray*} Defining $s=ir$ and $t=d/r$,
this becomes \begin{eqnarray*} & & \exp \left( - \sum_{s \geq 1} \sum_{t \geq
1} \frac{q^t u^{st}}{st q^{st}} \sum_{r|s} \mu(r) (-1)^{s/r}
(1-q^{-rt})^{-s/r} \right)\\
& = & \exp \left( - \sum_{s \geq 1} \sum_{t \geq
1} \frac{q^t u^{st}}{st q^{st}} \sum_{r|s} \mu(r) (-1)^{s/r} \sum_{b \geq 0} {s/r+b-1 \choose b} q^{-rtb} \right). \end{eqnarray*} Letting $a=rb$, this becomes \begin{eqnarray*}
&  & \exp \left( - \sum_{a \geq 0} \sum_{s \geq 1}
\sum_{t \geq 1} \frac{(u^sq^{1-s-a})^t}{t} \frac{1}{s} \sum_{r|s,a}
\mu(r) (-1)^{s/r} {s/r+a/r-1 \choose a/r} \right)\\ & = & \prod_{a \geq 0}
\prod_{s \geq 1} (1-u^sq^{1-s-a})^{F(a,s)}. \end{eqnarray*} \end{proof}

    As a consequence one has the following result, which is sharp in
    the sense that there are values of $d$ (such as $d=4$) for which
    the assertion would be false if the upper bound on $n$ does not
    hold. Parts 2 and 3 of Theorem \ref{stabcor} are known from
    \cite{W2}.

\begin{theorem} \label{stabcor}
\begin{enumerate}
\item The numbers \[F(a,s) := \frac{1}{s} \sum_{r|s,a} \mu(r)
(-1)^{s/r} {s/r + a/r - 1 \choose a/r} \] are integers for all $a \geq
0, s \geq 1$.
\item The coefficients of $q^{-n}$ in $s_{\GL(d,q)}$ and
$s_{\GL(\infty,q)}$ are equal whenever \mbox{$n \leq d-1$}.
\item The coefficients of $q^{-n}$ in $c_{\GL(d,q)}$ and
$c_{\GL(\infty,q)}$ are equal whenever \mbox{$n \leq 2d$}.
\item The coefficients of $q^{-n}$ in $s_{\U(d,q)}$ and $s_{\U}(\infty,q)$ are equal whenever \mbox{$n \leq d-1$}.
\item The coefficients of $q^{-n}$ in $c_{\U(d,q)}$ and $c_{\U(\infty,q)}$ are equal whenever \linebreak[4] \mbox{$n \leq 2d$}.
\end{enumerate}
\end{theorem}
\begin{proof} For
the first assertion, it follows from Theorem \ref{GeneralIntegers2}
that all coefficients $u^iq^{-j}$ in $(1+\frac{u}{q-1}) s_{\GL}(u,q)$
are integers. Now consider the expression for $(1+\frac{u}{q-1})
s_{\GL}(u,q)$ in Proposition \ref{stab1}. If some $F(a,s)$ were
non-integral, let $(a,s)$ be the smallest such, where smallest means
to first compare the $s$ coordinate, then if necessary the $a$
coordinate. Then the coefficient of $u^s q^{1-s-a}$ in
$(1+\frac{u}{q-1}) s_{\GL}(u,q)$ would be non-integral, a
contradiction.

    For the second assertion, note by M\"{o}bius inversion that
    $F(0,1)=-1$, $F(0,2)=1$ and that $F(0,m)=0$ for $m \geq 2$. Thus
    by Proposition \ref{stab1}, \[ (1-u) s_{\GL}(u,q) =
    \frac{(1-u^2q^{-1})}{(1+\frac{u}{q-1})} \prod_{a,s \geq 1} (1-u^s
    q^{1-s-a})^{F(a,s)}.\] The integrality of the $F(a,s)$ implies
    that $q^{-d}$ divides the coefficient of
    $u^{d+1}$ in $(1-u) s_{\GL}(u,q)$. This coefficient is
    $s_{\GL(d+1,q)}-s_{\GL(d,q)}$, which proves the result.

    For the third assertion, it is proved in \cite{W2} that \[
    c_{\GL(d+1,q)} - c_{\GL(d,q)} = q^{-d-1} \left[ s_{\GL(d+1,q)}
    - s_{\GL(d,q)} \right]. \] Thus by the previous paragraph, $q^{-2d-1}$
    divides $c_{\GL(d+1,q)} - c_{\GL(d,q)}$ as a polynomial, which implies the result.

    For the fourth assertion, note from Theorem 2.1.13 of \cite{FNP}
    that \[ s_\U(u,q) = \frac{s_{\GL(u^2,q^2)}}{s_{\GL(-u,-q)}}.\]
    From the expression for $(1-u) s_{\GL}(u,q)$ in the proof of the
    second assertion, it follows that $(1-u) s_\U(u,q)$ is equal to \[
    (1-u^2q^{-1}) \frac{\left( 1+\frac{u}{q+1} \right)}{\left(
    1+\frac{u^2}{q^2-1} \right)} \prod_{a,s \geq 1} \left(
    \frac{(1-u^{2s} q^{2(1-s-a)})}{(1+(-1)^{a}u^{s} q^{1-s-a})}
    \right)^{F(a,s)} .\] Thus Proposition \ref{stab1} implies that
    $q^{-d}$ divides the coefficient of $u^{d+1}$ in \mbox{$(1-u)
    s_{\U}(u,q)$}. This coefficient is $s_{\U(d+1,q)}-s_{\U(d,q)}$, as
    desired.

    For the fifth assertion, note from Theorem 2.1.10 of
    \cite{FNP} that \[ c_{\U(d+1,q)} - c_{\U(d,q)} = (-q)^{-d-1}
    \left[ s_{\U(d+1,q)} - s_{\U(d,q)} \right]. \] Thus by the previous
    paragraph, $q^{-2d-1}$ divides $c_{\U(d+1,q)} -
    c_{\U(d,q)}$, which implies the result. \end{proof}

\subsection{A General Stabilization Result} \label{general}

This subsection gives an approach to finding the rate
of stabilization of the finite dimensional coefficients to the
limiting coefficients which is more general, in that it is effective
for all $\Lambda$-types, though it does not give the sharpest possible
results in all cases. We describe this approach only in the case of
the groups $\GL(d,q)$, but it could be extended without difficulty to
$\G\in\{\U,\Symp,O,O^{\pm}\}$.

We shall need the following extension of Lemma \ref{qTransform1}:
\begin{lemma}\label{SubsetLemma}
Let $S$ be a subset of $\mathbb{N}\times\mathbb{N}$, closed under (componentwise) addition. 
Suppose that for integers $a_{i,j}$
$$
1+\sum_{1\le i,j} a_{i,j}u^{i}q^{-j}=\prod_{1\le i,j}(1-u^{i}q^{-j})^{b_{i,j}}.
$$
Then $a_{i,j}=0$ for all $(i,j)\notin S$ if and only if $b_{i,j}=0$ for all $(i,j)\notin S$.
\end{lemma}
\begin{proof}
It is clear that the product
$$
\prod_{\stack{1\le i,j}{(i,j)\in S}}(1-u^{i}q^{-j})^{b_{i,j}},
$$
when expressed as a power series in $u$ and $q^{-1}$, will only yield terms $u^{i}q^{-j}$ when 
$(i,j)\in S$. This is enough to prove one half of the double implication; the other half follows easily from induction on $n:=i+j$, using the explicit construction of $b_{k, n-k}$ given in the proof of Lemma \ref{qTransform1}. 
\end{proof}

Let $\Lambda_{0}$ be the set of all partitions of positive integers, and suppose that
\mbox{$\emptyset\neq \Lambda\subseteq\Lambda_{0}$}. We find lower bounds for the rate of stabilization of the coefficients of the polynomials $\Lambda_{\GL(d,q)}$ as $d$ increases. We use two similar methods, one for the case when $(1)\in\Lambda$, and the other for the case $(1)\notin \Lambda$. It is worth remarking that this particular distinction is intuitively reasonable; if $(1)\in \Lambda$, then all separable transformations are of $\Lambda$-type, and it follows that $\Lambda_{\GL(\infty,q)}\ge
s_{\GL}(\infty,q)>0.$ But if $(1)\notin\Lambda$, it is easy to show---indeed our argument will show---that $\Lambda_{\GL(\infty,q)}=0$. In the first case, our method will be to look at the difference $\Lambda_{\GL(d,q)}-\Lambda_{\GL(d-1,q)}$, and show that it is divisible (as a polynomial) by a particular power of $q^{-1}$. In the second case, we show that $\Lambda_{\GL(d,q)}$ itself is divisible by a power of $q^{-1}$.

For a non-empty partition $\lambda$, define $\Delta(\lambda)$ to be the degree of 
$C_{\GL,\lambda}(q)$ as a polynomial in $q$. 
This degree may be expressed in several ways:
\begin{lemma}\begin{enumerate}
\item Suppose that $\lambda$ has $m_{i}$ parts of size $i$ for all $i$. Then
$\Delta(\lambda)=2\sum_{i<j}im_{i}m_{j}+\sum_{i}im_{i}^{2}$.
\item Suppose that $n_{i}=\sum_{j\ge i} m_{j}$ for all $i$. Then
$\Delta(\lambda)=\sum_{i}n_{i}^{2}$.
\item Finally, suppose $\lambda$ has parts $a_{1},\dots, a_{k}$, where $a_{i}\ge a_{i+1}$ for all $i$. Then $\Delta(\lambda)=\sum_{i}(2i-1)a_{i}=2\sum_{i}ia_{i} - |\lambda|$ 
\end{enumerate}\end{lemma}
\begin{proof}
Definition \ref{centralizers} yields the first equation easily. The
second follows from the first, via the observation that
$\sum_{i<j}im_{i}m_{j}=\sum_{k\ge 1}\sum_{k\le i < j}m_{i}m_{j}$. The
third follows from the first by observing that, if $a_{s+1},\dots
a_{s+m_{i}}$ are the parts of size $i$, then
$\sum_{k=1}^{m_{i}}(2(s+k)-1)a_{s+k}=2ism_{i}+ im_{i}^{2}$.  Now use the
fact that $s=\sum_{j>i}m_{j}$, and sum over $i$.
\end{proof}

It can be established easily (using any of the expressions for $\Delta(\lambda)$ above) that if $\#(\lambda)$ denotes the number of parts of $\lambda$, then 
$$
|\lambda|\le\Delta(\lambda)\le|\lambda|\#(\lambda)\le|\lambda|^{2}.
$$
Each of these inequalities may in fact be equality, and in each case this imposes a regular structure on $\lambda$; in particular, we remark that 
$|\lambda|=\Delta(\lambda)$ if and only if $\lambda$ has a single part.

Define 
$$
T_{\Lambda}(u,q):=\sum_{\lambda\in\Lambda}\frac{u^{|\lambda|}}{C_{\GL,\lambda}(q)}.
$$
Then by the cycle index of $\GL(n,q)$ (\cite{St}, \cite{F})
\begin{equation}\label{Lambda1}
1+\sum_{d \geq 1} \Lambda_{\GL(d,q)}u^{d} =  
\prod_{d \geq 1 }\left(1+T_{\Lambda}(u^{d},q^{d})\right)^{N(d,q)}. 
\end{equation}
Setting $\Lambda=\Lambda_{0}$ shows that 
$$
\prod_{d \geq 1}\left(1+T_{\Lambda_{0}}(u^{d},q^{d})\right)^{N(d,q)} = \frac{1}{1-u}.
$$ 
If we write $\Lambda^{c}$ for the complement $\Lambda_{0}\setminus\Lambda$, then we obtain from (\ref{Lambda1}) the following equation:
\begin{equation}
\label{Lambda2} 
1+\sum_{d \geq 1} \left(\Lambda_{\GL(d,q)}-\Lambda_{\GL(d-1,q)}\right)u^{d}  =  
\prod_{d \geq 1} \left(
1-\frac{T_{\Lambda^{c}}(u^{d},q^{d})}{1+T_{\Lambda_{0}}(u^{d},q^{d})}\right)^{N(d,q)}.
\end{equation} Here $\Lambda_{\GL(0,q)}$ is to be interpreted as 1.

Now $1+T_{\Lambda_{0}}(u,q)$ may be written in the form $1+\sum_{1\le
i\le j} r_{i,j}u^{i}q^{-j}$. Its reciprocal can also be put into this
form, since the modification of Lemma \ref{qTransform1} in which all
occurrences of $1 \leq i,j$ are replaced by $1 \leq i \leq j$ is
true. In fact it is shown in \cite{St} that
$\left(1+T_{\Lambda_{0}}(u,q)\right)^{-1} = \prod_{r \geq 1}
(1-uq^{-r})$.  It follows that there are integers $a_{i,j}$ such that
$$
1-\frac{T_{\Lambda^{c}}(u,q)}{1+T_{\Lambda_{0}}(u,q)}=1+\sum_{2 \le
i,j} a_{i,j}u^{i}q^{-j}.
$$ 
Note the condition $2\le i,j$ on the index of summation, which is valid since 
$\Lambda^c$ does not contain the partition $(1)$. 
In fact it is easily shown that if $a_{i,j}\neq 0$, then
there exists a partition $\lambda$ in $\Lambda^{c}$ such that $i\ge
|\lambda|$, $j\ge \Delta(\lambda)$, and $j-i\ge
\Delta(\lambda)-|\lambda|$. Let $S$ be the set of all pairs $(i,j)$
satisfying this condition. Then $S$ is obviously closed under
addition, and it follows from Lemma \ref{SubsetLemma} that we can find
integers $b_{i,j}$ such that
$$
1-\frac{T_{\Lambda^{c}}(u,q)}{1+T_{\Lambda_{0}}(u,q)} = \prod_{(i,j)\in S}
(1-u^{i}q^{-j})^{b_{i,j}}.
$$
Now by our usual argument, invoking Lemma \ref{Identities}, part \ref{GL-},
$$
\prod_{d \geq 1} \left(1-\frac{T_{\Lambda^{c}}(u^{d},q^{d})}
{1+T_{\Lambda_{0}}(u^{d},q^{d})}\right)^{N(d,q)}=\prod_{(i,j)\in S}\left(\frac{1-u^{i}q^{1-j}}{1-u^{i}q^{-j}}\right)^{b_{i,j}},
$$
which may certainly be put into the form
$$
\prod_{\stack{i,j \geq 1} {(i,j+1)\in S}}(1-u^{i}q^{-j})^{c_{i,j}}
$$
for integers $c_{i,j}$. 

Define $\sigma:= \mathrm{inf}\{\,j/i\,\mid \,(i,j+1)\in S\}$. Then the set
$\{(i,j)\,\mid \,j\ge i\sigma\}$ is closed under addition. By Lemma \ref{SubsetLemma} it follows that, for some integers $e_{i,j}$, we may write
$$ 
\prod_{d \geq 1} \left(1-\frac{T_{\Lambda^{c}}(u^{d},q^{d})}
{1+T_{\Lambda_{0}}(u^{d},q^{d})}\right)^{N(d,q)}=1+\sum_{\stack{i,j \geq 1}{j\ge i\sigma}}e_{i,j}u^{i}q^{-j}.
$$
Then from (\ref{Lambda2}) above, it follows that for all $d$,
$$
\Lambda_{\GL(d,q)}-\Lambda_{\GL(d-1,q)}=\sum_{j\ge d\sigma}e_{d,j}q^{-j},
$$ 
and hence that, whenever $j<(d+1)\sigma$, the coefficient of $q^{-j}$ in $\Lambda_{\GL(d,q)}$ has stabilized to the coefficient in the limit $\Lambda_{\GL(\infty,q)}$.

 What is the ratio $\sigma$? Firstly, suppose that $\Lambda^{c}$ contains at least one partition with a single part, and suppose that the smallest such partition is $(k)$. Then $\Delta\left(\,(k)\,\right)=k$, and it is clear (since $\Delta(\lambda)>|\lambda|$ for partitions with more than one part) that $\sigma=\frac{k-1}{k}$. Suppose, on the other hand, that $\Lambda^{c}$ contains no one-part partition. Then it is clear that 
$\frac{j-1}{i}\ge 1$ for any $(i,j)\in S$. But it is also clear that by taking $i$ and $j$ sufficiently large, we can make this ratio arbitrarily close to $1$, and hence that $\sigma=1$.      

We summarize these conclusions in the following theorem:
\begin{theorem}\label{stabthe1}
Suppose $\{\,(1)\,\}\subseteq \Lambda\subseteq \Lambda_{0}$. Define
$$
\sigma:=\left\{\begin{array}{cl}
1-\frac{1}{k} & \textrm{if $k$ is the size of the smallest one-part partition not in $\Lambda$,} \\
1 & \textrm{if $\Lambda$ contains all one-part partitions.}
\end{array}\right.
$$
Then the coefficient of $q^{-j}$ in 
$\Lambda_{\GL(d,q)}$ is equal to the coefficient in $\Lambda_{\GL(\infty,q)}$ whenever 
$j<(d+1)\sigma$.
\end{theorem}

We have the following applications to cyclic, separable and semisimple matrices:
\begin{cor}
\begin{enumerate}
\item The coefficients of $q^{-j}$ in $c_{\GL(d,q)}$ and $c_{\GL(\infty,q)}$
are equal whenever $j\le d$.
\item The coefficients of $q^{-j}$ in $s_{\GL(d,q)}$ and $s_{\GL(\infty,q)}$
are equal whenever $j\le \frac{d}{2}$.
\item The coefficients of $q^{-j}$ in $ss_{\GL(d,q)}$ and $ss_{\GL(\infty,q)}$
are equal whenever $j\le \frac{d}{2}$.
\end{enumerate} \end{cor}
Comparing with Theorem \ref{stabcor}, these bounds are not sharp for
the cases of $s_{\GL(d,q)}$ and $c_{\GL(d,q)}$. However it is not at
all clear that the methods of \cite{LeS}, \cite{W2}, or Theorem \ref{stabcor} can be adapted to
the semisimple (or other) cases.

We now deal with the case when $(1)\notin\Lambda$. Let $k$ be the smallest integer such that $\Lambda$ contains a partition of $k$. For each $i\ge k$, define 
$$
\tau_{i}:=\mathrm{min}\{\frac{\Delta(\lambda)}{|\lambda|}\,\mid\,\lambda\in\Lambda, \, 
|\lambda|\le i\}.
$$
Write
$$
1+T_{\Lambda}(u,q)=1+\sum_{2 \le i,j}a_{i,j}u^{i}q^{-j}.
$$
Note that $i,j \geq 2$ since $(1) \notin \Lambda$. Also it is not hard to see that $a_{i,j}=0$ unless $j\ge i\tau_{i}$. Furthermore, since 
$\tau_{i}$ is weakly decreasing function of $i$, the set $S:=\{(i,j)\,\mid\,j\ge i\tau_{i}, \ i,j \geq 2\}$ is 
additively closed. By Lemma \ref{SubsetLemma}, there are integers $b_{i,j}$ such that
$$
1+T_{\Lambda}(u,q)=\prod_{(i,j)\in S}(1-u^{i}q^{-j})^{b_{i,j}}.
$$
Now, proceeding as usual by means of part (\ref{GL-}) of Lemma \ref{Identities}, we find that
$$
\prod_{d \geq 1}
\left(1+T_{\Lambda}(u^{d},q^{d})\right)^{N(d,q)}=\prod_{(i,j)\in S}
\left(\frac{1-u^{i}q^{1-j}}{1-u^{i}q^{-j}}\right)^{b_{i,j}},
$$
which can be rewritten in the form
\begin{equation}\label{cform}
\prod_{(i,j+1)\in S}(1-u^{i}q^{-j})^{c_{i,j}}.
\end{equation}
Suppose $c_{i,j}\neq 0$. Then there is a partition $\lambda_{i}$ such that 
$|\lambda_{i}|\le i$, and $\frac{j+1}{i}\ge\frac{\Delta(\lambda_{i})}{|\lambda_{i}|}$.
But now
$$
\frac{j}{i}\ge \frac{\Delta(\lambda_{i})}{|\lambda_{i}|}-\frac{1}{i}
\ge \frac{\Delta(\lambda_{i})}{|\lambda_{i}|}-\frac{1}{|\lambda_{i}|}
=\frac{\Delta(\lambda_{i})-1}{|\lambda_{i}|}\ge \sigma,
$$
where $\sigma:=\mathrm{inf}\{\frac{\Delta(\lambda)-1}{|\lambda|}\,\mid\,\lambda\in\Lambda\}$.
The set $\{(i,j)\,\mid\,j\ge i\sigma, \  i,j \geq 1 \}$ is additively closed, and so (\ref{cform}) may be written, by Lemma \ref{SubsetLemma}, in the form
$$
\sum_{\stack{i,j \geq 1}{j\ge i\sigma}}e_{i,j}u^{i}q^{-j}, 
$$
for integers $e_{i,j}$. 

Now from (\ref{Lambda1}) above, the expression at (\ref{cform}) is equal to $1+\sum_{d}\Lambda_{\GL(d,q)}u^{d}$. It follows that 
$$
\Lambda_{\GL(d,q)}=\sum_{j\ge d\sigma}e_{d,j}q^{-j},
$$
which suffices to prove the following theorem:
\begin{theorem}\label{stabthe2}
Suppose that $\emptyset\neq\Lambda\subseteq\Lambda_{0}$ and that $(1)\notin\Lambda$. Define
$$
\sigma:=\mathrm{inf}\{\frac{\Delta(\lambda)-1}{|\lambda|}\,\mid\,\lambda\in\Lambda\}.
$$
Then whenever $j<d\sigma$, the coefficient of $q^{-j}$ in $\Lambda_{\GL(d,q)}$ is $0$.
\end{theorem}

The constant $\sigma$ is likely to be fairly easy to calculate for most naturally arising sets $\Lambda$. If $\Lambda$ contains one-part partitions, and $(k)$ is the smallest, then $\sigma=1-\frac{1}{k}$. If $\Lambda$ has no one-part partitions, then $\sigma\ge 1$.

\end{document}